\documentclass[12pt]{article}  

\usepackage{amsmath}
\usepackage{amssymb}
\usepackage{cite}
\usepackage{paralist}
\usepackage[usenames]{color}
\usepackage{graphicx}          
\usepackage[dvips]{epsfig}    
\newtheorem{definition}{Definition}
\newtheorem{theorem}{Theorem}

\title{Partial stabilization of stochastic systems with application to rotating rigid bodies\thanks{
This work was supported in part by the State Fund for Fundamental Research of Ukraine (project F78/206-2018) and the budget program of NAS of Ukraine (KPKVK 6541230).\newline
$^{1}$Max Planck Institute for Dynamics of Complex Technical Systems, Magdeburg, Germany
\newline
$^{2}$Institute of Applied Mathematics and Mechanics, National Academy of Sciences of Ukraine
{\tt\small (e-mail: zuyev@mpi-magdeburg.mpg.de, irisna.shurko@gmail.com)}
}
}

\author{Alexander Zuyev$^{1,2}$ and Iryna Vasylieva$^{2}$}
\date{}

\begin{document}

\maketitle
\thispagestyle{empty}

\begin{abstract}
This paper addresses the problem of stabilizing a part of variables for control systems described by stochastic differential equations of the It\^o type.
The considered problem is related to the asymptotic stability property of invariant sets and has important applications in mechanics and engineering.
Sufficient conditions for the asymptotic stability of an invariant set are proposed by using a stochastic version of LaSalle's invariance principle.
These conditions are applied for constructing the state feedback controllers in the problem of single-axis stabilization of a rigid body. The cases of control torques generated by jet engines and rotors are considered as illustrations of the proposed control design methodology.
\end{abstract}

\section{Introduction}

The concept of stability with respect to a part of variables was introduced by Lyapunov and systematically investigated in the monographs by~\cite{RO87},~\cite{Vo98},~\cite{Zuyev_15}, and other authors.

An effective method for the study of stability problems for systems with random actions is the stochastic Lyapunov's direct method.
In the paper by~\cite {Sh_78}, stability conditions with respect to a part of variables and conditions for  partial asymptotic stability of It\^o's stochastic differential equations have been obtained.
The work by \cite {Soch_18} addresses stability issues for hybrid systems described by stochastic differential equations with three types of Markovian switching processes. For these processes, criteria of exponential $y$-stability have been obtained in the linear and nonlinear framework.

The paper by \cite{Ig_13} is devoted to the study of partial asymptotic stability for stochastic differential equations with time-invariant vector fields.
A modification of the Barbashin--Krasovskii theorem for the case of stochastic asymptotic stability of the trivial equilibrium with respect to a part of variables has been proposed in the above-mentioned paper.
The concepts of stability of solutions with respect to a part of variables and stability of an invariant  set, being very similar, are not identical.
For systems with random actions, namely, for a system of the It\^o-type differential equations, \cite {Stan_01} obtained conditions for the invariance of sets and their stochastic stability. Note that these results do not imply the asymptotic stability of an invariant set, and the study of limit sets for solutions to the stochastic differential equations remains a challenging problem (see, e.g.,~\cite{Oks_03} and references therein).
For dynamical systems with multivalued flows on a metric space, asymptotic stability conditions of invariant sets have been derived in the paper by~\cite{Zu03}.

Our goal is to propose asymptotic stability conditions for invariant sets of nonlinear stochastic differential equations
and to apply such conditions to the partial stabilization problem of stochastic control systems.
This paper is organized as follows.
The notions of partial stability and partial asymptotic stability are introduced in Section~2.
Sufficient conditions for partial asymptotic stability are presented in Section~3.
These results are applied for the single-axis stabilization of a rotating rigid body with jet controls and rotors (flywheels) in Sections~4 and~5, respectively.

\section{Partial stability of stochastic systems}

In this section, we present some basic definitions related to partial stability of stochastic systems.

Consider a system of stochastic differential equations of the It\^o type:
\begin{equation}
dx(t)=f(x)dt+\sigma(x)dW(t),\quad x\in {\mathbb R}^n.
\label{sys1}
\end{equation}
Here $x=(x_1,...,x_n)^T$ is the state vector of the system, the functions    $f:\mathbb{R}^n \to \mathbb{R}^n$ and  $\sigma:\mathbb{R}^n \to \mathbb{R}^{n\times k}$ are assumed to be  Lipschitz continuous on each compact set $K\subset \mathbb{R}^n$. We treat $X=(x_1,...,x_n)^T$ and $f(x)=(f_1(x),...,f_n(x))^T$ as column vectors, and $\sigma(x)$ as $n\times k$ matrix.
System~\eqref{sys1} is subject to the $k$-dimensional Wiener process $W(t)$ whose components $w_j(t)$ $(j = 1,2, ..., k)$ are independent one-dimensional Wiener processes.
Under these conditions, there exists a unique strictly Markov process $x(t; x^0,s)$ which is a solution of system~\eqref{sys1}  under the initial condition $x(s; x^0, s)=x^0$.

\begin{definition}\label{def1}
A nonempty set $M \subset \mathbb{R}^n$ is called:
\begin{itemize}
\item {\em forward invariant} for system~\eqref{sys1} if $x^0 \in M $ and $s \in \mathbb{R}$ imply $x(t; x^0, s)$ for all $t\ge s$ almost surely;
\item {\em backward invariant} for system~\eqref{sys1} if $x^0 \in M $ and $s \in \mathbb{R}$ imply $x(t; x^0, s)$ for all $t\le s$ almost surely;
\item {\em invariant} for system~\eqref{sys1} if it is forward and backward invariant.
\end{itemize}
\end{definition}

We assume in the sequel that the state vector of system~\eqref{sys1} can be written as $x=(y^T,z^T)^T$ with  $y=(y_1,...,y_m)^T \in \mathbb{R}^m$ and  $z=(z_1,...,z_p)^T\in \mathbb{R}^p$, $m+p=n$.

\begin{definition}\label{def2}
The set $M=\{x\in{\mathbb R}^n\,|\,y=0\}$ is called stable in probability if, for any $s \ge 0,$ $\epsilon_1 >0$ and $\epsilon_2 >0$, there exists a $\delta >0$ such that the following inequality holds for all initial data $x^0 \in \mathbb{R}^n$ with $||y^0||<\delta:$
\begin{equation}
P\{\sup_{t\ge s} ||y(t;x^0,s)||>\epsilon_1\}<\epsilon_2.
\label{stab_property}
\end{equation}
\end{definition}

\begin{definition}\label{def3}
The set $M=\{x\in{\mathbb R}^n\,|\,y=0\}$ is called asymptotically stable in probability if it is stable in probability and, for any $\varepsilon>0$, there is a $\Delta=\Delta(\varepsilon)>0$ such that
\begin{equation}
P\{\lim_{t \to +\infty}||y(t;y^0 ,s)||=0\}>1-\varepsilon,\;\forall x^0\in{\mathbb R}^n:\, \|y^0\|<\Delta.
\end{equation}
\end{definition}

In order to study the stability of invariant sets, we consider the class of functions $V:\mathbb{R}^n\rightarrow\mathbb{R}^{+} =[0;+\infty),$ $V(0)=0$, which are twice continuously differentiable.
We relate with system~\eqref{sys1} the differential operator
\begin{equation}\mathcal L=\sum_{i=1}^n f_{i}(x)\frac{\partial }{\partial x_{i}}+\frac{1}{2}\sum_{i,j=1}^{n}a_{ij}(x)\frac{\partial^{2}}{\partial x_{i}\partial x_{j}},\end{equation} where $[a_{ij}]=\sigma\sigma^{T}$.
Let us also introduce the class $\mathcal {K}$ of comparison functions, consisting of all continuous strictly increasing functions $\alpha:\mathbb{R}^{+}\rightarrow\mathbb{R}^{+}$ such that $\alpha(0)=0$.

\section{Main result}

In this section, sufficient conditions for the asymptotic stability of invariant sets of system ~\eqref{sys1} will be derived.

\begin{theorem}\label{thm1}
Let the  set $M=\{x\in {\mathbb R}^n\,|\,y=0\}$ be invariant for system~\eqref{sys1},
and let $V\in \mathcal C^2(\mathbb R^{n}; \mathbb R^+)$ satisfy the following properties:
\begin{itemize}
\item[1)] $\alpha_{1}(\|y\|)\leq V(x)\leq \alpha_{2}(\|y\|)$ for all $x \in \mathbb R^{n}$ with some $\alpha_1, \alpha_2 \in \mathcal K$;
\item[2)]$\mathcal LV(x)\leq 0$ for all $x \in \mathbb R^{n}$;
\item[3)] there exists a $\Delta > 0$ such that $x(t;x^0,s)$ is bounded for $t\ge s$ almost surely, provided that $||y^0||< \Delta$;
\item[4)] The  set $\{x\in {\mathbb R}^n\,|\, \mathcal LV(x)=0\} \backslash M$ is not forward invariant for~\eqref{sys1}.
\end{itemize}

Then the set $M=\{x\in {\mathbb R}^n\,|\,y=0\}$ is asymptotically stable in probability for system~\eqref{sys1}.
\end{theorem}

{\em Proof.}
Let us take an arbitrary $\varepsilon >0.$ Let $||y^0||<\varepsilon $, and let $$\tau_{\varepsilon} = \inf_{t \ge s}\{t\,|\, ||y(t; x^0, s)|| > \varepsilon\},\; \tau_{\varepsilon}(t)= \min(\tau_{\varepsilon},t).$$
Using a special case of  Dynkin's formula (see, e.g.,~\cite{Oks_03}), we obtain:
$${\mathbb E} V(x(\tau_{\varepsilon}(t);x^0,s)) = V(x^0) + {\mathbb E}\int\limits_s^{\tau_{\varepsilon}(t)} \mathcal LV(x(u;x^0,s))du.
$$
Therefore,
\begin{equation}{\mathbb E}  V(x(\tau_{\varepsilon}(t);x^0,s)) \le V(x^0)\quad \text{for all}\; t>s,\label{Eineq}
\end{equation}
where ${\mathbb E}$ is the expectation with respect to the probability measure $P_{x^0,s}\{\cdot\}$

Let us rewrite inequality~\eqref{Eineq} in the following form:
$$\int\limits_{\tau_{\varepsilon}<t} \alpha_1 (||y(\tau_{\varepsilon};x^0 ,s)||)P_{x^0,s}(d\omega)+$$$$+\int\limits_{\tau_{\varepsilon}\ge t} \alpha_1 (||y(\tau_{\varepsilon};x^0 ,s)||)P_{x^0,s}(d\omega)\le V(x^0).$$

So, $\alpha_1(\varepsilon)P_{x^0,s}\{\tau_{\varepsilon}<t\}\le V(x^0).$
This implies, because of  the continuity of $V (x)$ and the equality $V (0)=0$, that
$$\lim_{y^0 \to 0}P_{x^0,s}\{\tau_{\varepsilon}<t\}=0,$$
so,~\eqref{stab_property} holds.

It remains to prove that $$\lim_{y^0 \to 0}P\{\lim_{t \to \infty}||y(t;x^0 ,s)||=0\}=1.$$
For this purpose we select the subset $\mathcal B$ of the set of sample trajectories of the process $x(t; x^0,s)$ such that $\tau_{\varepsilon}(t)=t$
holds for each component $y_i(t; x^0,s)$ of
$x(t; x^0,s)$  for all $t\ge s$.

Then
\begin{equation} \lim_{y^0\rightarrow 0}P(\mathcal B)=1.\end{equation}

By the assumption~3) of Theorem~1, there exists a $ \Delta > 0$ such that all solutions $x(t;x^0,s)$ of system~\eqref{sys1} with $||y^0||< \Delta$ are bounded for $t>s$.
Let $||y^0||< \Delta,$ then let us show that $\lim_{t \rightarrow \infty}||y(t; x^0 ,s)||=0.$

Due to the boundedness of $x(t;x^0,s)$, there exists a sequence $t_k \rightarrow +\infty:$
$$ x(t_k; x^0,s)\rightarrow x^{*}=(y^{*},z^{*}) \in G^{+,}$$ where $G^{+}$ is the set of $\omega$-limit  points of the solution $x(t; x^0,s).$

Since the conditions of the stochastic version of LaSalle's invariance principle are satisfied (see~\cite{Mao_99,Mao_00}), $G^{+}\subset \{x \,|\,\mathcal LV(x)=0\}$.

Now we show that $y^{*} = 0$. By the invariance of the set $G^{+}$, we obtain $x(t;x^{*},s) \in G^{+}$  for all $t \ge s$. Suppose that $y^{*} \ne 0$. Since the set $M$ is invariant and $x^{*}\in G^{+}\subset \{x\,|\, \mathcal LV(x)=0\}$, then $x(t;x^{*},s) \in \{x\,|\,  \mathcal LV(x)=0\} \backslash M$, $t \ge s$.

The above property contradicts the assumption~4) of Theorem~1.
Thus, $y^{*} = 0$ and \begin{equation}\lim_{t \rightarrow +\infty}||y(t; x^0 ,s)||=0.
\label{y_limit}
\end{equation}

Then~\eqref{y_limit} implies the assertion of Theorem~1.
$\square$

\section{Single-axis  stabilization of a satellite by jet torques}

Consider the problem of stabilizing the  motion of a rigid body around its center of mass under the action of jet control torques without taking into account the change of mass. The equations of motion can be written in the Euler--Poisson form with random disturbances as follows:
\begin{equation}
\begin{array}{lcl}
d\omega_{1}=(\frac{A_{2}-A_{3}}{A_{1}}\omega_2\omega_3+u_1)dt+\omega_1\sigma
dW(t),\\
d\omega_{2}=(\frac{A_{3}-A_{1}}{A_{2}}\omega_1\omega_3+u_2)dt+\omega_2\sigma
dW(t),\\
d\omega_{3}=(\frac{A_{1}-A_{2}}{A_{3}}\omega_1\omega_2)dt,\\ d\nu_{1}=(\omega_3\nu_2-\omega_2\nu_3)dt,\\
d\nu_{2}=(\omega_1\nu_3-\omega_3\nu_1)dt,\\
d\nu_{3}=(\omega_2\nu_1-\omega_1\nu_2)dt.
\end{array}
\label{Euler1}
\end{equation}
Here $(\omega_1,\omega_2,\omega_3  )$ are projections of the  angular velocity vector $\omega$ on the corresponding principal axes of inertia of the body, $(\nu_1,\nu_2,\nu_3  )$  are projections of a fixed vector $\nu$ onto the corresponding principal inertia axes, $(A_1, A_2, A_3 )$ are  principal central moments of inertia of the body, and $(u_1,u_2  )$ are control torques. In the further analysis, we will assume that all variables and parameters are dimensionless.

System~\eqref{Euler1} with $u_1=u_2=0$ admits the following particular solution:
\begin{equation}\begin{array}{c}\omega_1=\omega_2=\omega_3=0,\\
\nu_1=\nu_2=0, \nu_3=1.\end{array}
\label{sol1}
\end{equation}

Solution~\eqref{sol1} corresponds to the equilibrium position at which the third principal axis of inertia of the body is directed along the $\nu$ vector.  We apply Theorem~1 to stabilize the set $M=\{(\omega_1,\omega_2,\omega_3,\nu_1,\nu_2,\nu_3)|\nu_1=\nu_2=\omega_1=\omega_2=0\}$ of system~\eqref{Euler1}. This goal  corresponds to the problem of single-axis stabilization of the rigid body, i.e. the projections $\nu_1,\nu_2$ and their derivatives  $\dot\nu_1,\dot\nu_2$ should be small and tending to zero as $t\rightarrow+\infty,$ while the other components of the solutions are assumed to be bounded.

Let us take the following control Lyapunov function candidate:$$V=\frac{1}{2}(A_1\omega_1^2+A_2\omega_2^2)+\frac{1}{2}(\nu_1^2+\nu_2^2).$$

Then we define the feedback control as follows:\begin{equation}\begin{aligned}
u_1=\omega_2\omega_3-\frac{1}{A_{1}}\nu_2\nu_3-\Big\{\frac{|A_1-A_2|}{2A_1}|\omega_3|+A_1h+\frac{1}{2}\sigma\sigma^{T}\Big\}\omega_1,
\\u_2=-\omega_1\omega_3+\frac{1}{A_{2}}\nu_1\nu_3-\Big\{\frac{|A_1-A_2|}{2A_2}|\omega_3|+A_2h+\frac{1}{2}\sigma\sigma^{T}\Big\}\omega_2,\end{aligned}
\label{controls1}
\end{equation}
where the function $h$ will be defined below.

By computing the action of the $\mathcal L$ operator on $V$, we get: $$\mathcal LV=-\frac{|A_1-A_2|}{2}|\omega_3|\omega_1^2-\frac{|A_1-A_2|}{2}|\omega_3|\omega_2^2-A_1^2h\omega_1^2-A_2^2h\omega_2^2\leq0.$$
Thus, the second condition of  Theorem 1 is satisfied.
All the solutions of system \eqref{Euler1} are bounded with respect to the variables $\nu_i$ due to the geometric integral $\nu_1^2+\nu_2^2+\nu_3^2 = const.$ It remains to verify the boundedness of solutions with respect to $\omega_i.$ For this purpose we exploit Theorem~39.1 from~\cite{RO87}.
According to that theorem,
it suffices to construct a function $W(x)$, which is positive definite with respect to  all   $\omega_i$, does not increase along the trajectories of the closed-loop system~\eqref{Euler1},~\eqref{controls1}, and satisfies the condition $W(x)\rightarrow\infty$ as $|\omega|\rightarrow\infty.$

Let us define $$W(x)=\frac{1}{2}(A_1\omega_1^2+A_2\omega_2^2+A_3\omega_3^2)+\frac{1}{2}(\nu_1^2+\nu_2^2).$$
Then $$\mathcal LW(x)=(A_1-A_2)\omega_1\omega_2\omega_3-\frac{|A_1-A_2|}{2}|\omega_3|(\omega_1^2+\omega_2^2)$$$$-h(A_1^2\omega_1^2+A_2^2\omega_2^2).$$
According to  Theorem~39.1 of~\cite{RO87}, for the $\omega_i$-boundedness of the solutions it suffices to satisfy the inequality $\mathcal LW_0(x)\leq0,$ so  the function $h(x)>0$  is chosen from the following condition:
\begin{equation}h(x)(A_1^2\omega_1^2+A_2^2\omega_2^2)\geq(A_1-A_2)\omega_1\omega_2\omega_3.
\label{ineq_h}
\end{equation}

Using the inequality $2A_1A_2\omega_1\omega_2\leq A_1^2\omega_1^2+A_2^2\omega_2^2,$  we conclude that to satisfy~\eqref{ineq_h} it suffices to take \begin{equation}h(x)=\left|\frac{A_1-A_2}{2A_1A_2}\omega_3\right|+\varepsilon,\end{equation} where $\varepsilon $ is any positive number.

Let us now check the invariance of the set $M$.
Thus we analyze the solutions of the system~\eqref{Euler1},~\eqref{controls1} for an arbitrary initial data from the set $M$. Let us rewrite the system \eqref{Euler1}, \eqref{controls1}, substituting $\nu_1=\nu_2=\omega_1=\omega_2=0:$
\begin{equation}
\begin{array}{lcl}
d\omega_{3}=0,\; d\nu_{3}=0.
\end{array}
\end{equation}
Then the solutions of the closed-loop system satisfy
\begin{equation}
\omega_{3}(t)=const\;\;\text{and}\;\;
\nu_{3}(t)=const
\end{equation}
almost surely.
 It remains to show that the condition~4) of Theorem~1 is satisfied. The set $M_v=\{x|\mathcal L V(x)=0\}$ has the form $M_v=\{x|\omega_1=\omega_2=0\}$, i.e. \begin{equation}x(t)\in M_v:
\begin{array}{lcl}
d\nu_1=\omega_3\nu_2dt,\;
d\nu_2=-\omega_3\nu_1dt,\\
0=-\frac{1}{A_1}\nu_2\nu_3dt,\\
0=\frac{1}{A_2}\nu_1\nu_3,\; d\omega_3=0.
\end{array}
\end{equation}

It is easy to verify that, for the initial conditions close enough to (9), all the trajectories of the closed-loop system~\eqref{Euler1} with~\eqref{controls1} satisfy the property $\nu_1=\nu_2=0.$
Moreover, the solutions of system (15) are $\nu_1(t)=0,\quad\nu_2(t)=0,\quad\omega_{3}(t)=const,\quad \nu_3(t)=const.$ The same solutions satisfy system (13).  So, $M_v \backslash M$ does not contain any positive trajectory of the considered closed-loop system.

All the conditions of Theorem 1 are satisfied, so the invariant set $M=\{(\omega_1,\omega_2,\omega_3,\nu_1,\nu_2,\nu_3)|\nu_1=\nu_2=\omega_1=\omega_2=0\}$ of the closed-loop system~\eqref{Euler1},~\eqref{controls1} is asymptotically stable in probability.

Below we present numerical simulations of the closed-loop system \eqref{Euler1}, \eqref{controls1} with the following dimensionless parameters: $A_1=1, A_2=3, A_3=2, \varepsilon = 0.01$.
These simulations have been performed in Maple by using the ItoProcess function.

\begin{figure}
\center {
\includegraphics[scale=0.3]{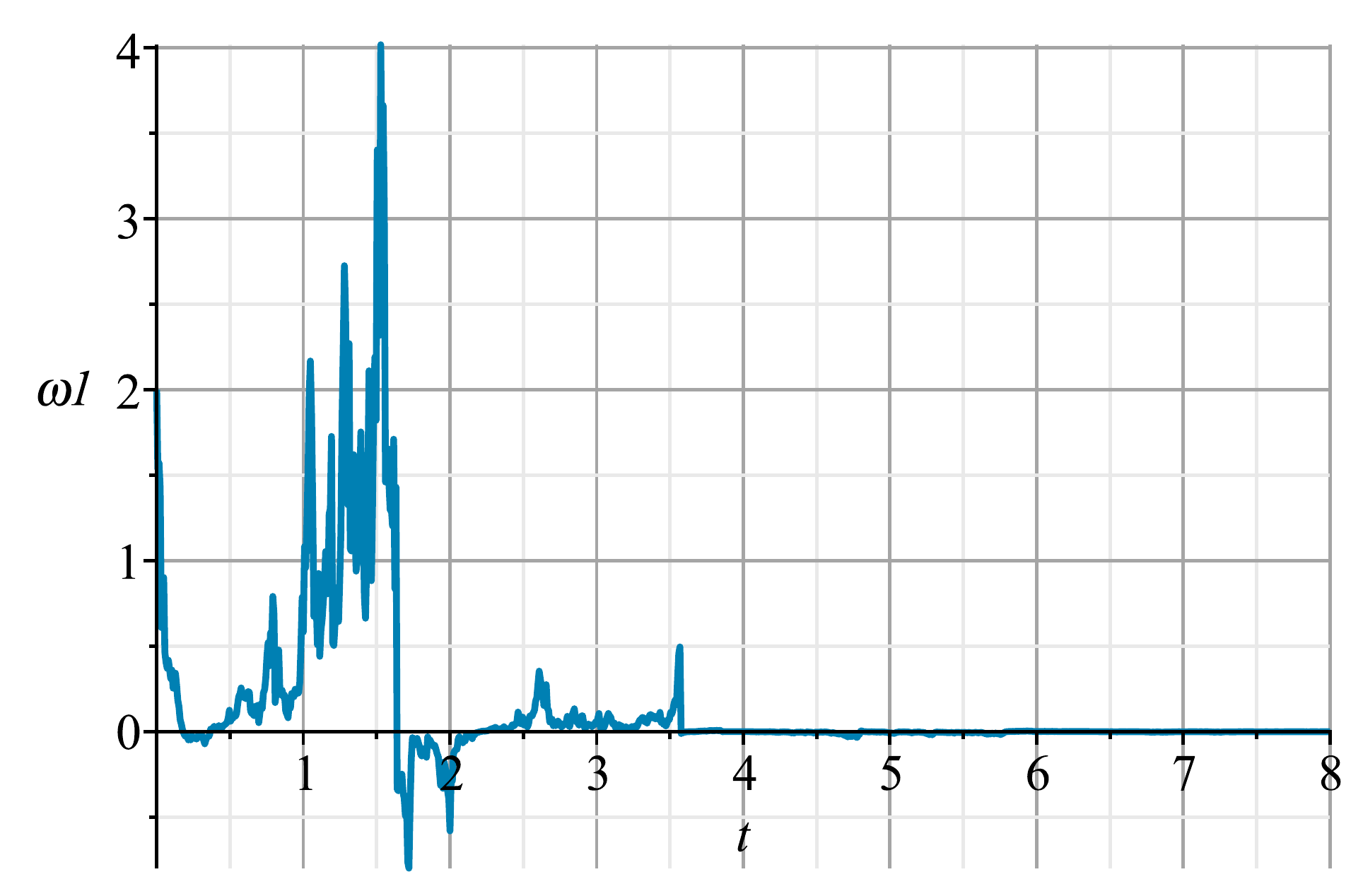}
\includegraphics[scale=0.3]{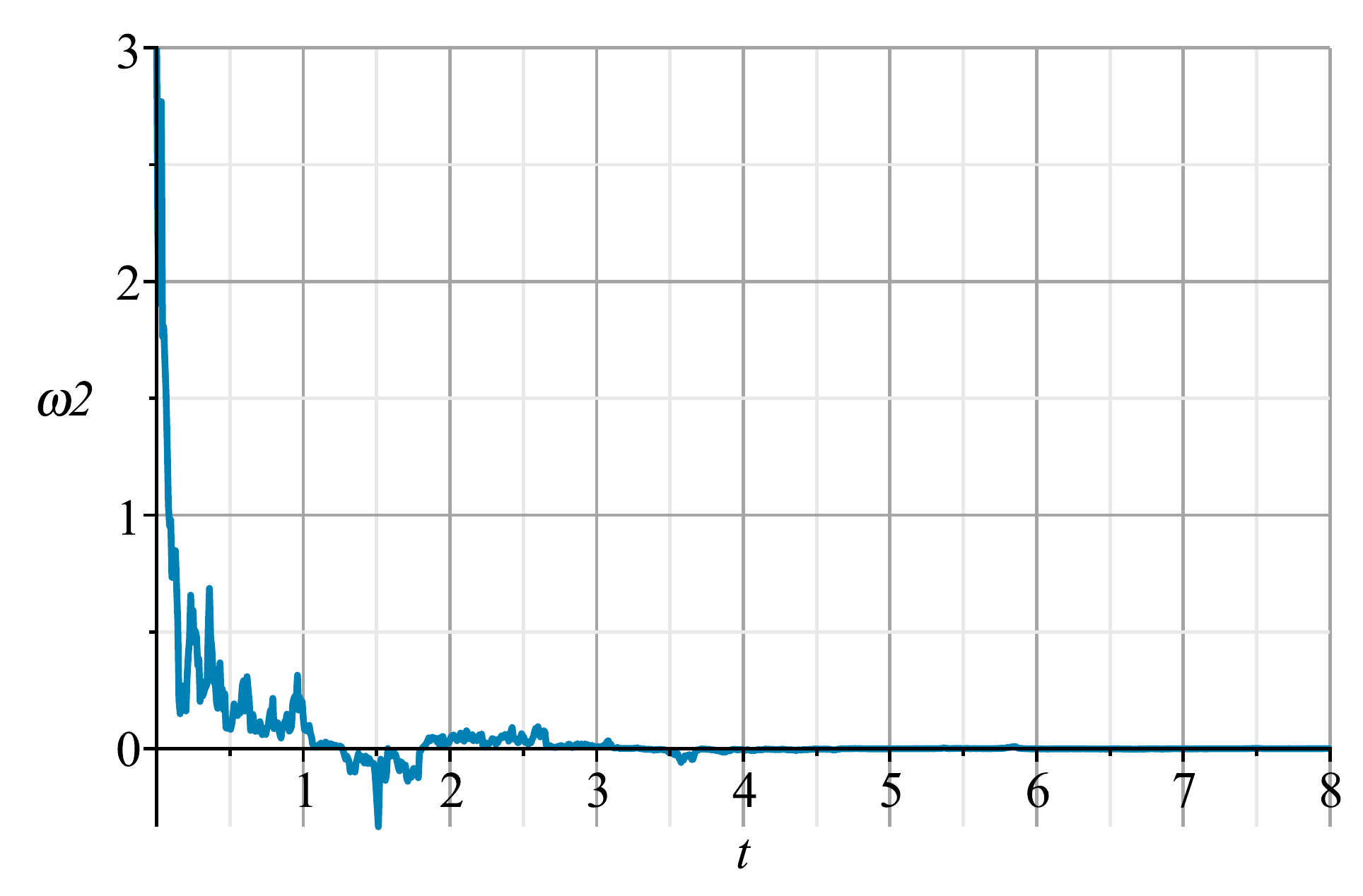}
\includegraphics[scale=0.3]{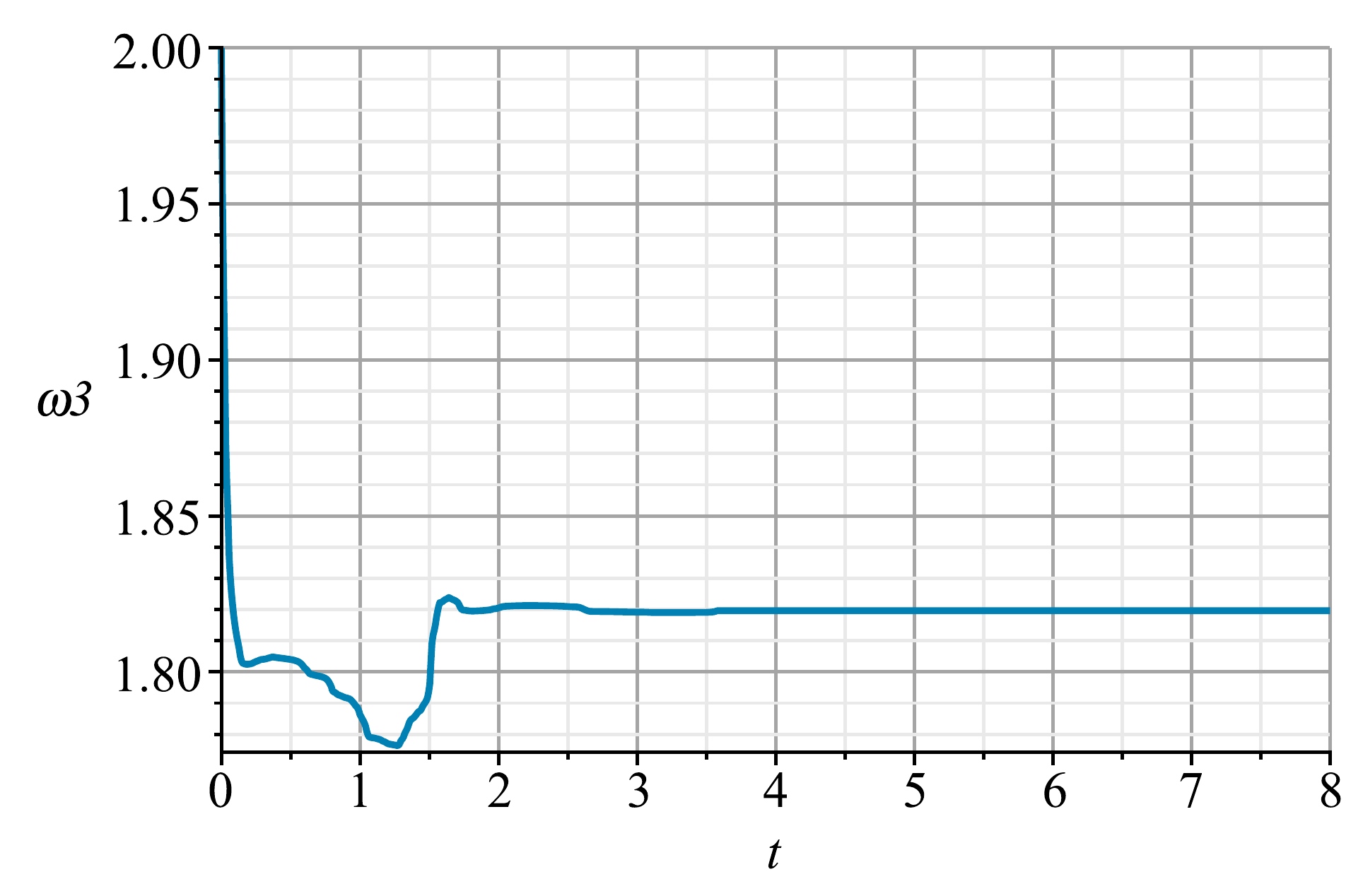}
}
\caption{Components $\omega_i(t)$ of a sample path of the closed-loop system~\eqref{Euler1}, \eqref{controls1}.}
\label{fig1w}
\end{figure}

\begin{figure}
\center {
\includegraphics[scale=0.3]{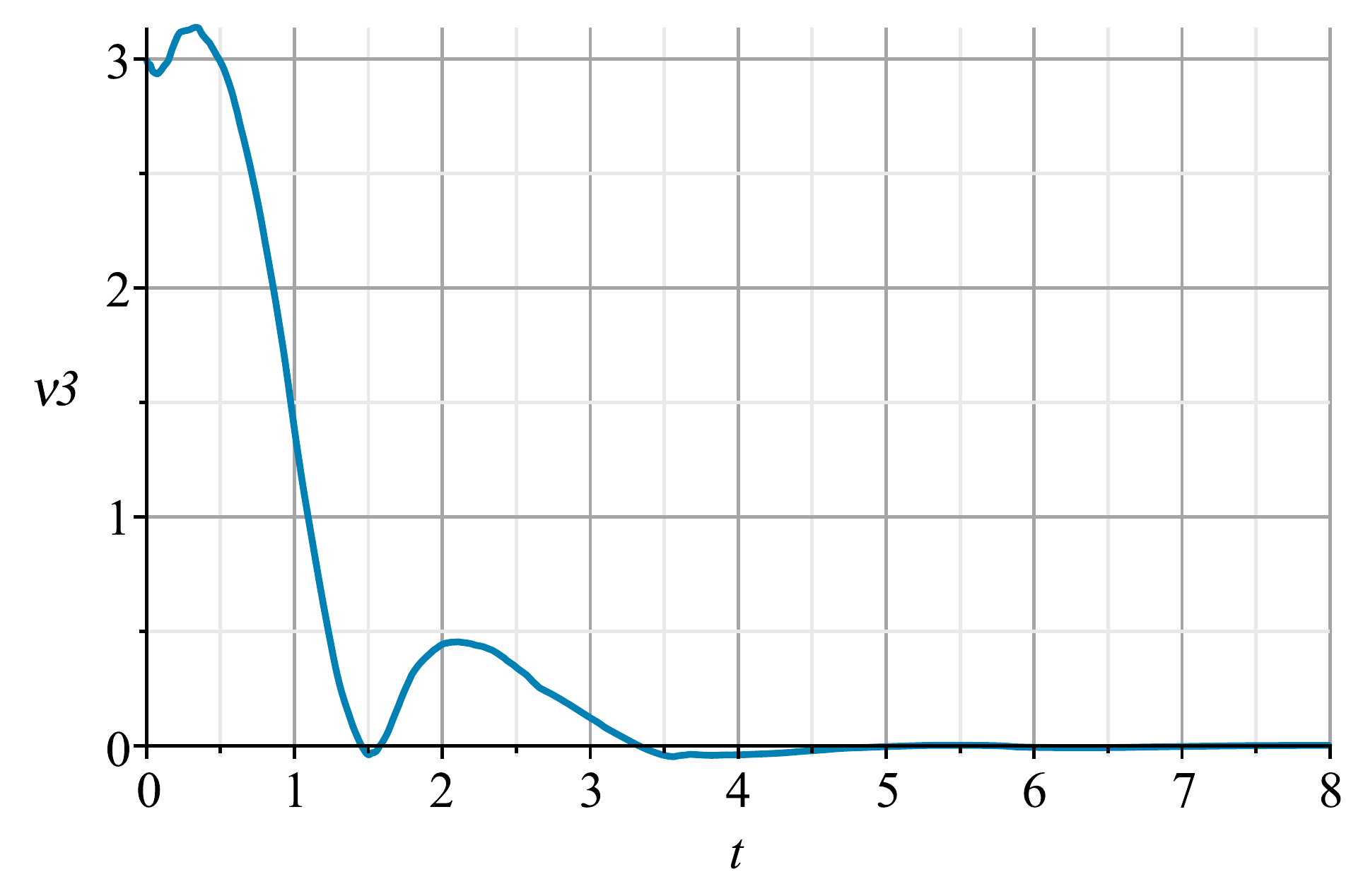}
\includegraphics[scale=0.3]{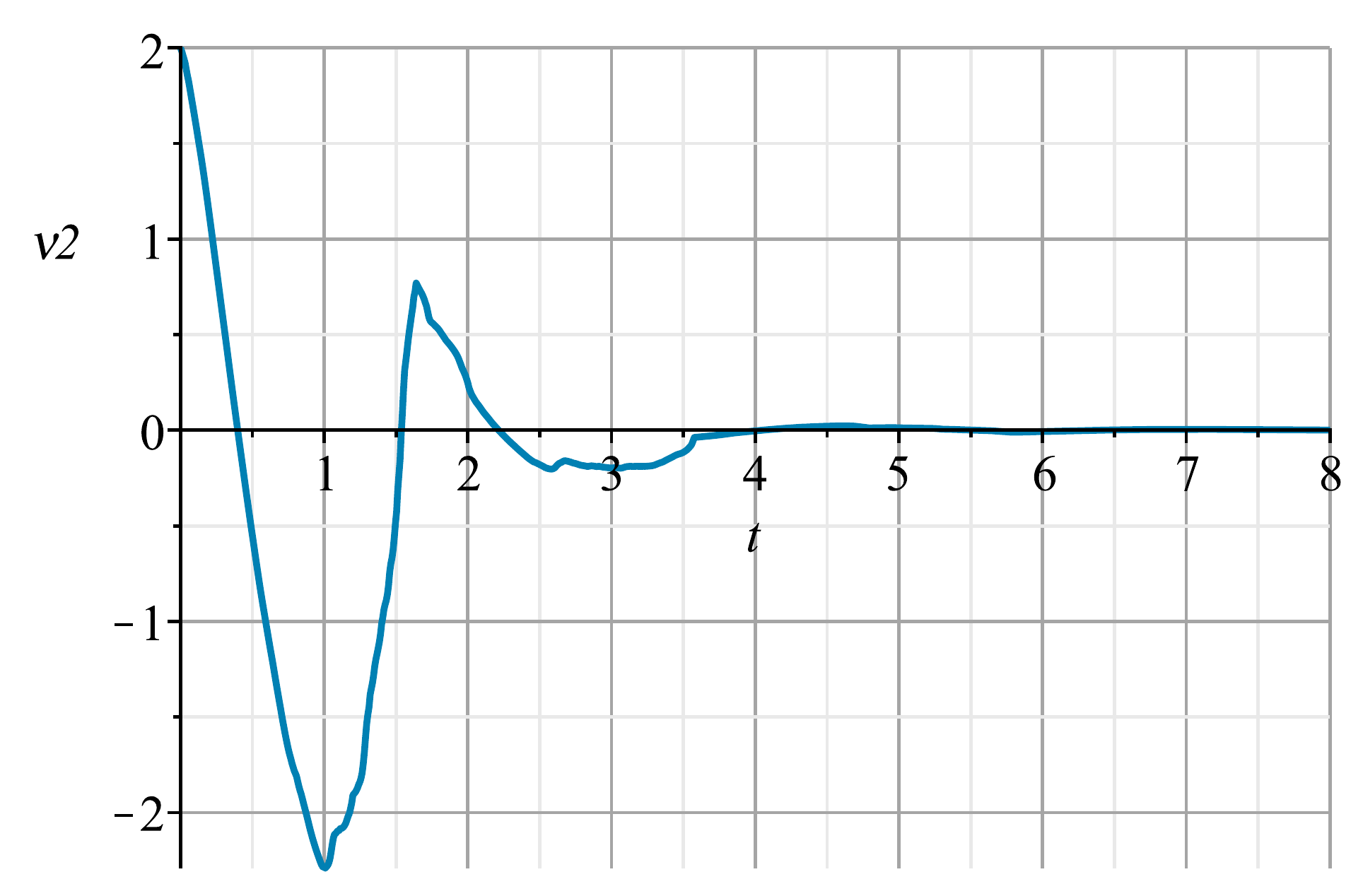}
\includegraphics[scale=0.3]{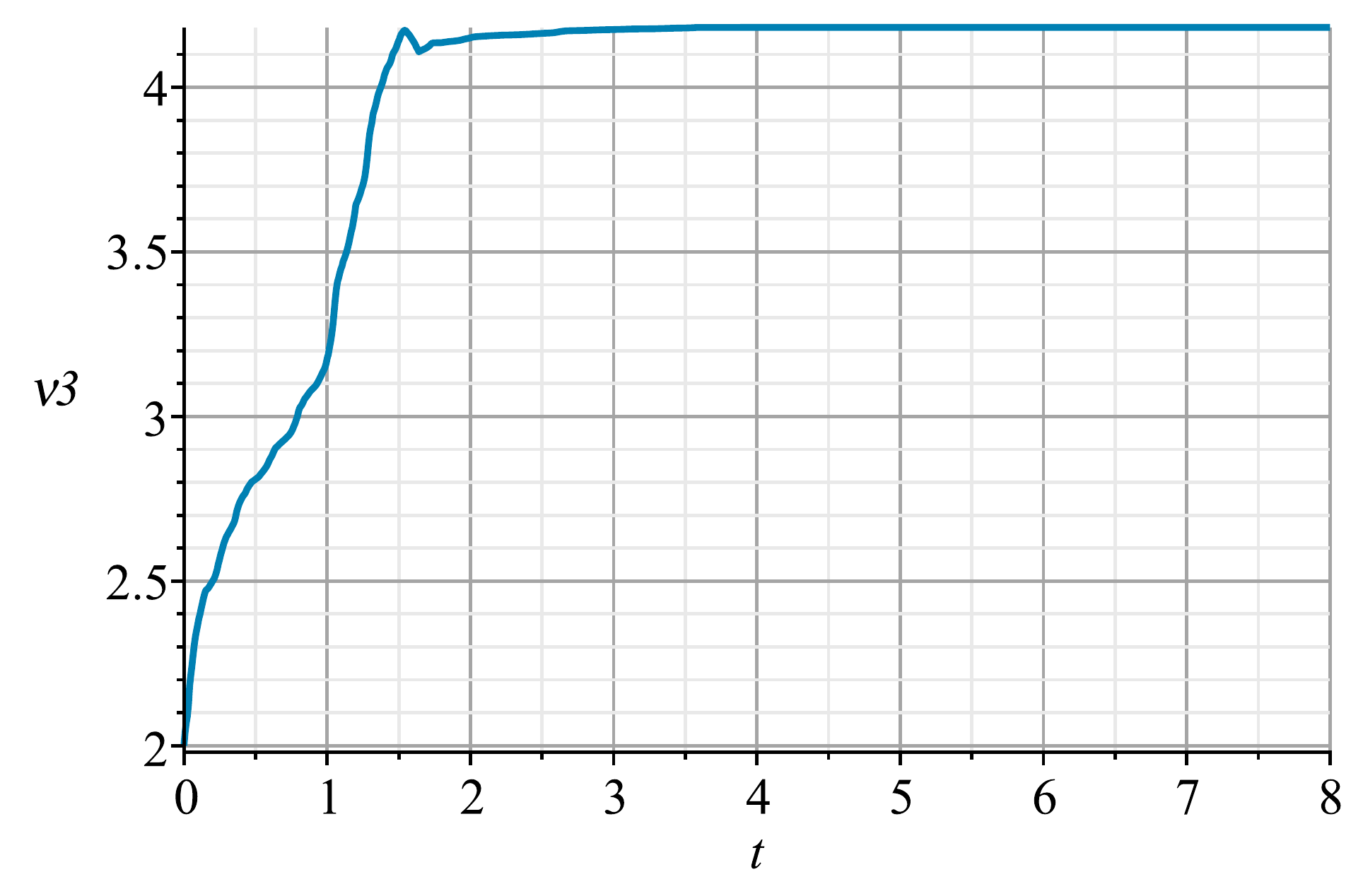}
}
\caption{Components $\nu_i(t)$ of a sample path of the closed-loop system~\eqref{Euler1}, \eqref{controls1}.}
\label{fig1nu}
\end{figure}
From Figs.~\ref{fig1w} and~\ref{fig1nu} it is clear that the components $\omega_1, \omega_2, \nu_1, \nu_2$ tend to zero for large values of $t$, while the coordinates $\omega_3, \nu_3$ tend to some limit values.
Thus, the limit motions of the body are uniform rotations around a fixed orientation vector, which is collinear to the third principal axis of inertia.

\section{Single-axis  stabilization of a satellite using two rotors}

The equations of motion of a rigid body containing a pair of symmetrical rotors with random effects can be written as follows:
\small\begin{equation}
\begin{array}{lcl}
d\omega_{1}=\left(\frac{A_2-A_3}{A_1-I_1}\omega_2\omega_3+\frac{I_2\Omega_2}{A_1-I_1}\omega_3-\frac{1}{A_1-I_1}u_1\right)dt +\omega_1\sigma dW(t),\\
d\omega_{2}=\left(\frac{A_3-A_1}{A_2-I_2}\omega_1\omega_3-\frac{I_1\Omega_1}{A_2-I_2}\omega_3-\frac{1}{A_2-I_2}u_2\right)dt+ \omega_2\sigma dW(t),\\
d\omega_{3}=\left(\frac{A_1-A_2}{A_3}\omega_1\omega_2+\frac{I_1\Omega_1}{A_3}\omega_2-\frac{I_2\Omega_2}{A_3}\omega_1\right)dt,\\
d\Omega_{1}=\left(\frac{1}{I_2}u_1-\frac{A_2-A_3}{A_1-I_1}\omega_2\omega_3-\frac{I_2\Omega_2}{A_1-I_1}\omega_3+\frac{1}{A_1-I_1}u_1\right)dt-\omega_1\sigma dW(t),\\
d\Omega_{2}=\left(\frac{1}{I_1}u_2-\frac{A_3-A_1}{A_2-I_2}\omega_1\omega_3+\frac{I_1\Omega_1}{A_2-I_2}\omega_3+\frac{1}{A_2-I_2}u_2\right)dt-\omega_2\sigma dW(t),\\
d\nu_{1}=(\omega_3\nu_2-\omega_2\nu_3)dt,\\
d\nu_{2}=(\omega_1\nu_3-\omega_3\nu_1)dt,\\
d\nu_{3}=(\omega_2\nu_1-\omega_1\nu_2)dt.
\end{array}
\label{Euler2}
\end{equation}\normalsize
Here $\omega_i  $ are coordinates of the angular velocity vector of the carrier body in the main coordinate frame, $\Omega_1,\Omega_2  $ are  relative angular velocities of the first and second rotor, respectively, $A_i $ are principal moments of inertia of the whole system consisting of the carrier body and rotors, $I_1,I_2  $ are moments of inertia of the first and second rotor, respectively, and $u_1,u_2$ are the control torques applied to the first and the second rotor, respectively.
 We assume that $A_1>I_1,A_2>I_2$.
System~\eqref{Euler2} is a stochastic version of the equations studied in~\cite{Zuyev_01}.
Note that the rotating rigid body with a rotor has been considered in the book on nonholonomic mechanics by~\cite{Bloch} as a mathematical model of a satellite.

Control system~\eqref{Euler2} with $u_1=u_2=0$ admits the following equilibrium: \begin{equation}\begin{array}{c}\omega_i=0, \quad \Omega_1=const, \quad \Omega_2=const,\\ \nu_1=\nu_2=0, \quad \nu_3=1.\end{array}\end{equation}

The considered mechanical system has the following integrals:\small
$$W_1 = (A_1\omega_1+I_1\Omega_1)^2+(A_2\omega_2+I_2\Omega_2)^2+(A_3\omega_3)^2=const,$$
$$W_2 = (A_1\omega_1+I_1\Omega_1)\nu_1+(A_2\omega_2+I_2\Omega_2)\nu_2+(A_3\omega_3)\nu_3=const,$$
$$W_3 =\nu_1^2+\nu_2^2+\nu_3^2=const.$$\normalsize

For the existence of the integral of moments, it is necessary for the vector $\sigma,$ which characterizes random actions, to satisfy the condition $\sigma \in N,$  where $N  $ is the invariant subspace of the linear operator $Q$ corresponding to zero eigenvalue. The linear operator $Q$ is given by the Hessian matrix of $W_1$:
\begin{equation}Q=\begin{pmatrix}
  2A_1^2 & 0 & 0 & 2A_1I_1 &0\\
  0& 2A_2^2 &0 &0&2A_2I_2\\
  0 &0 &2A_3^2 &0 &0\\
  2I_1A_1 &0 &0 &2I_2^2 &0\\
  0 &2I_2A_2 &0 &0 &2I_2^2
\end{pmatrix}.\end{equation}

Let us compute the  eigenvectors of $Q $ corresponding to zero eigenvalue. These vectors are
$$x=\begin{pmatrix}-\frac{I_1}{A_1}\Omega_1,-\frac{I_2}{A_2}\Omega_2,0,\Omega_1,\Omega_2\end{pmatrix}^T.$$

Thus the set  $N$ has the form $N=\{(\omega_1,\omega_2,\omega_3,\Omega_1,\Omega_2)^T|$ $$\omega_1=-\frac{I_1}{A_1}\Omega_1, \omega_2=-\frac{I_2}{A_2}\Omega_2, \omega_3=0\}.$$

To stabilize the set $M=\{(\omega_1,\omega_2,\omega_3,\Omega_1,\Omega_2,\nu_1,\nu_2,\nu_3)|\nu_1=\nu_2=\omega_1=\omega_2=0\}$ of system \eqref{Euler2}, we apply  Theorem 1 with the following control Lyapunov function candidate: $2V(x) = (A_1-I_1)\omega_1^2+(A_2-I_2)\omega_2^2+\nu_1^2+\nu_2^2.$

Let us define the feedback law as follows:
\begin{equation}
\begin{aligned}u_1=&\nu_2\nu_3+(A_2\omega_2+I_2\Omega_2)\omega_3\\&+\omega_1(h+\frac{1}{2}\sigma\sigma^T(A_1-I_1)+\frac{|A_1-A_2|}{2}|\omega_3|),\\
u_2=&-\nu_1\nu_3-(A_1\omega_1+I_1\Omega_1)\omega_3\\&+\omega_2(h+\frac{1}{2}\sigma\sigma^T(A_2-I_2)+\frac{|A_1-A_2|}{2}|\omega_3|).\end{aligned}
\label{controls2}
\end{equation}

\begin{figure}
\center{
\includegraphics[scale=0.3]{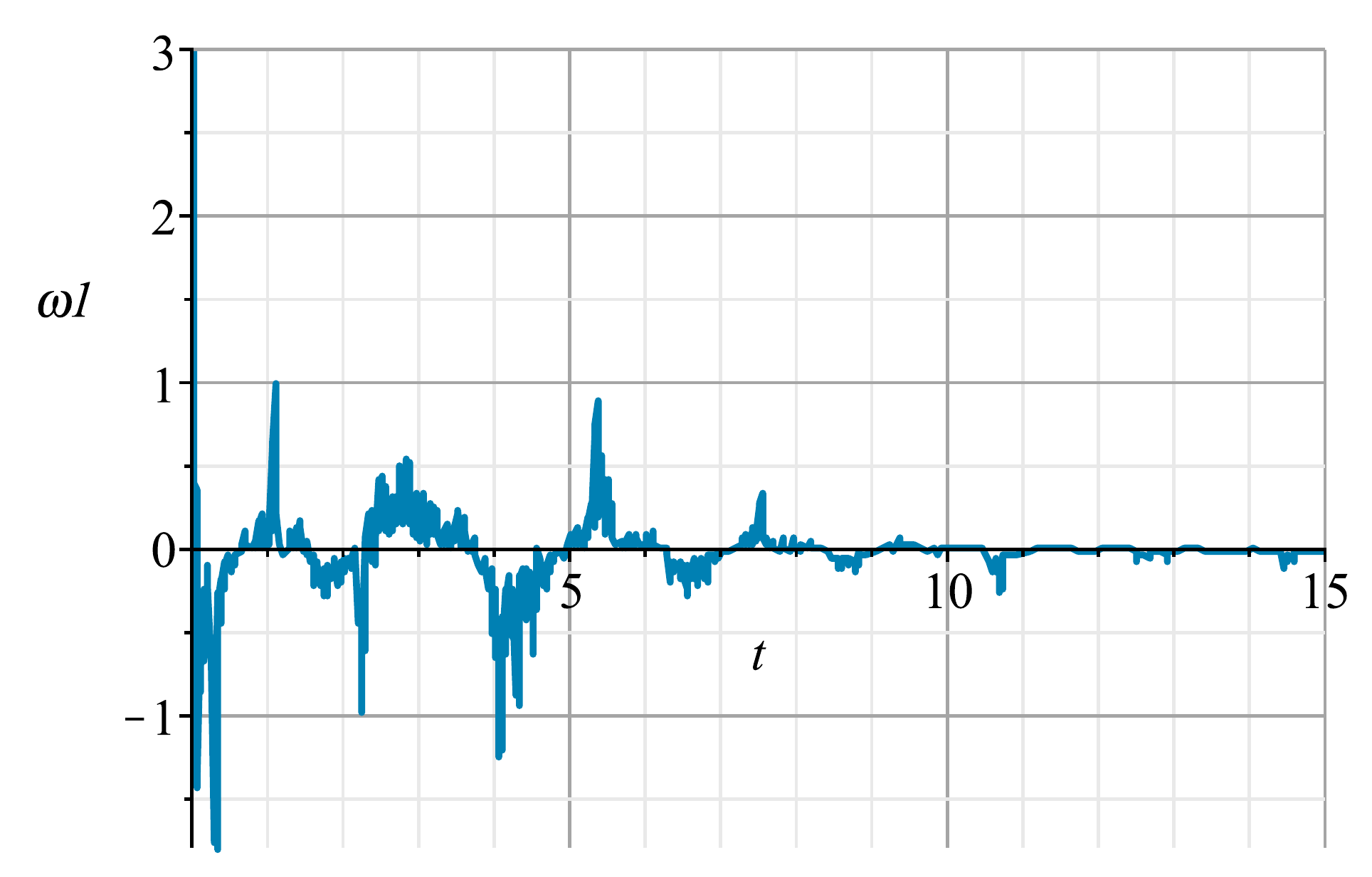}
\includegraphics[scale=0.3]{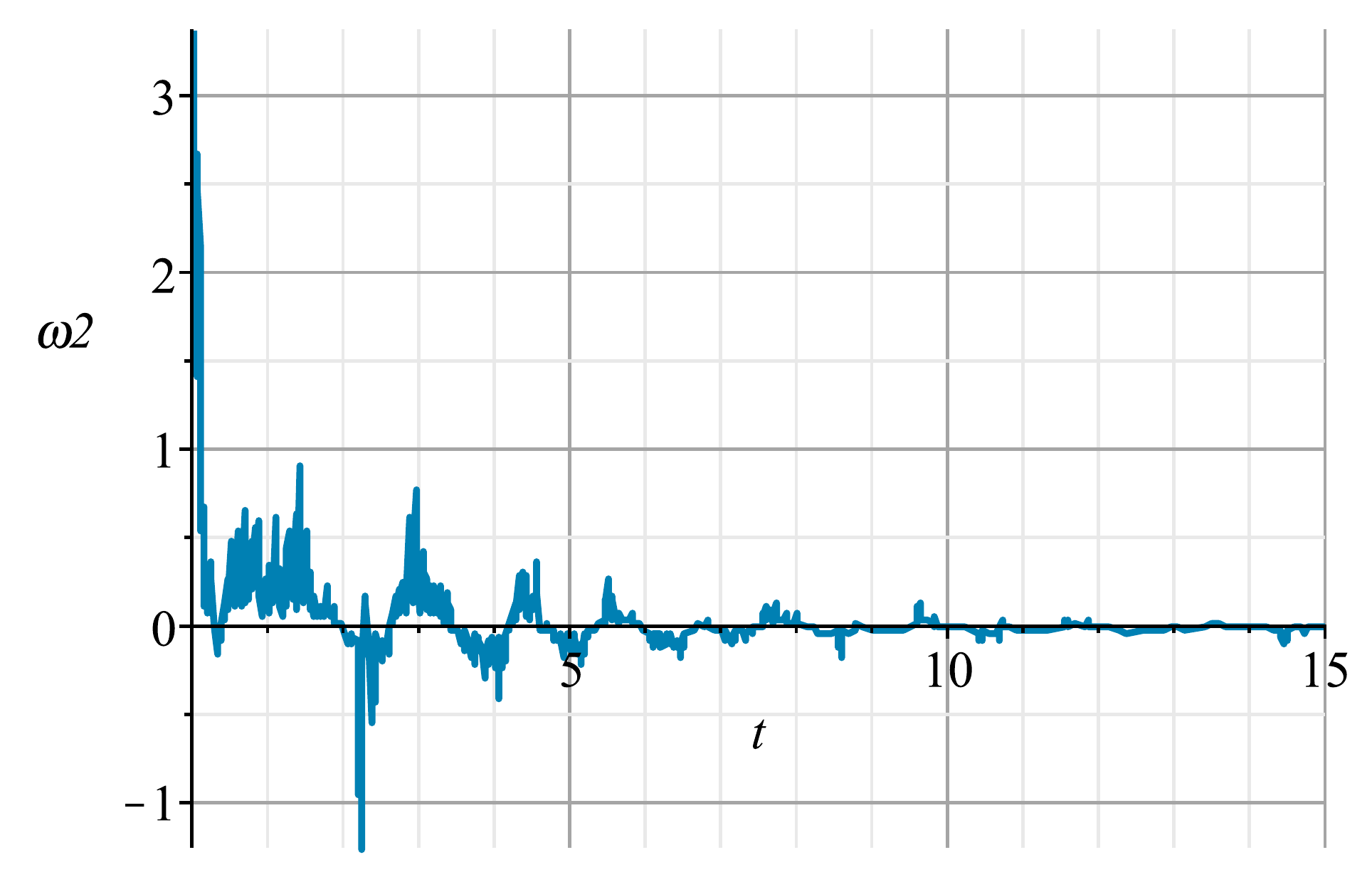}
\includegraphics[scale=0.3]{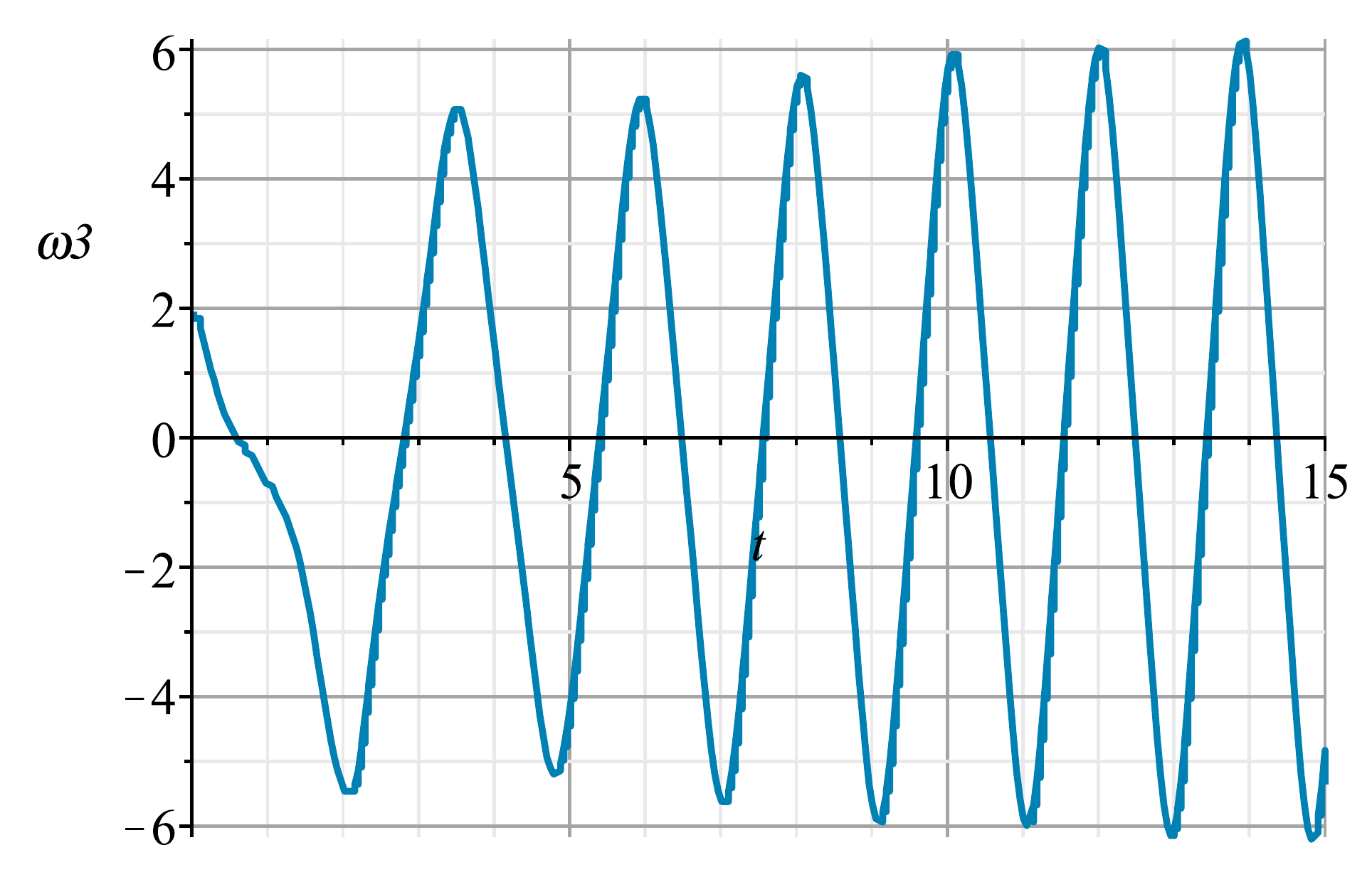}
}
\caption{Components $\omega_i(t)$ of a sample path of the closed-loop system~\eqref{Euler2}, \eqref{controls2}.}
\label{fig2w}
\end{figure}

\begin{figure}
\center{
\includegraphics[scale=0.3]{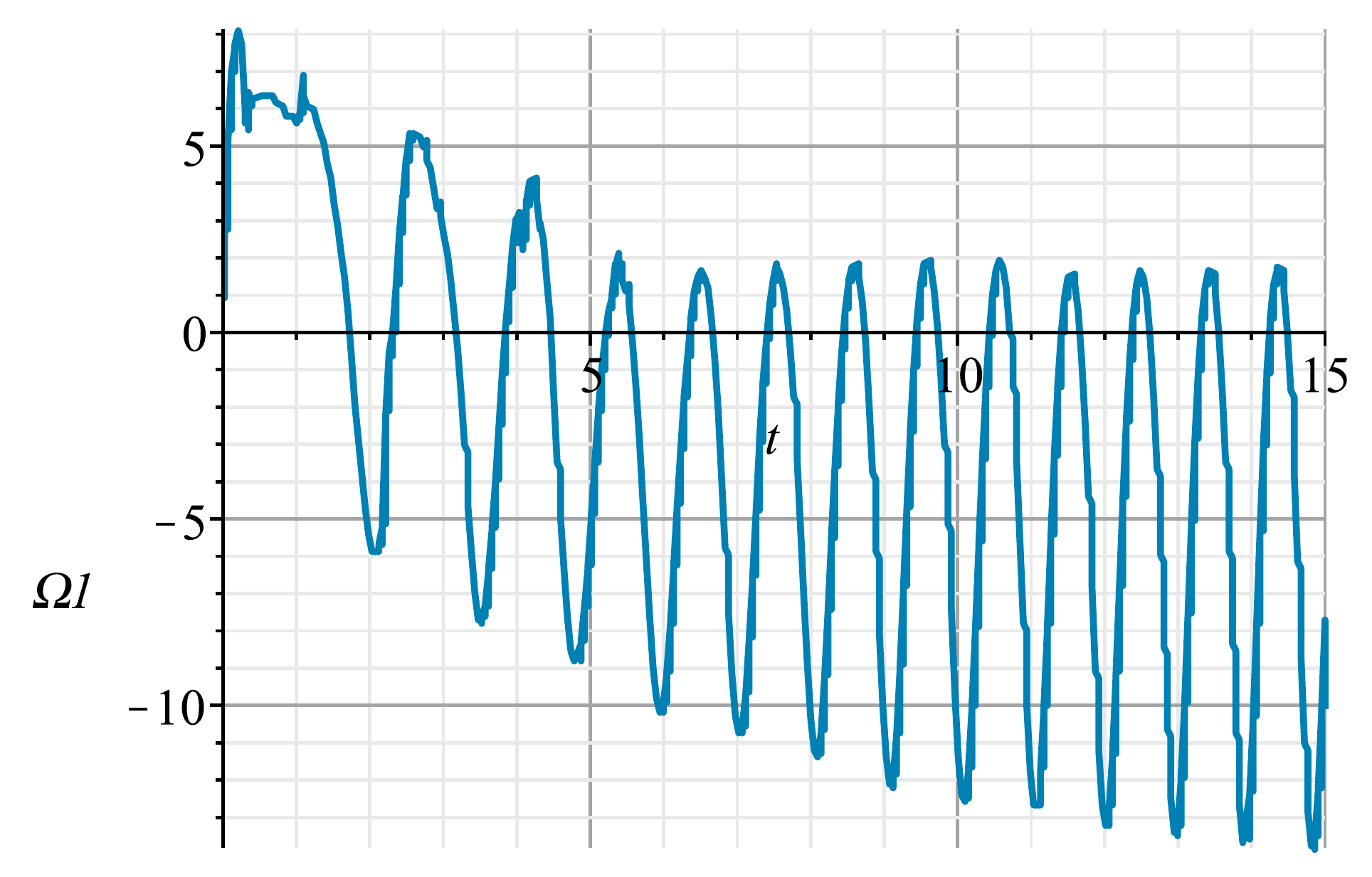}
\includegraphics[scale=0.3]{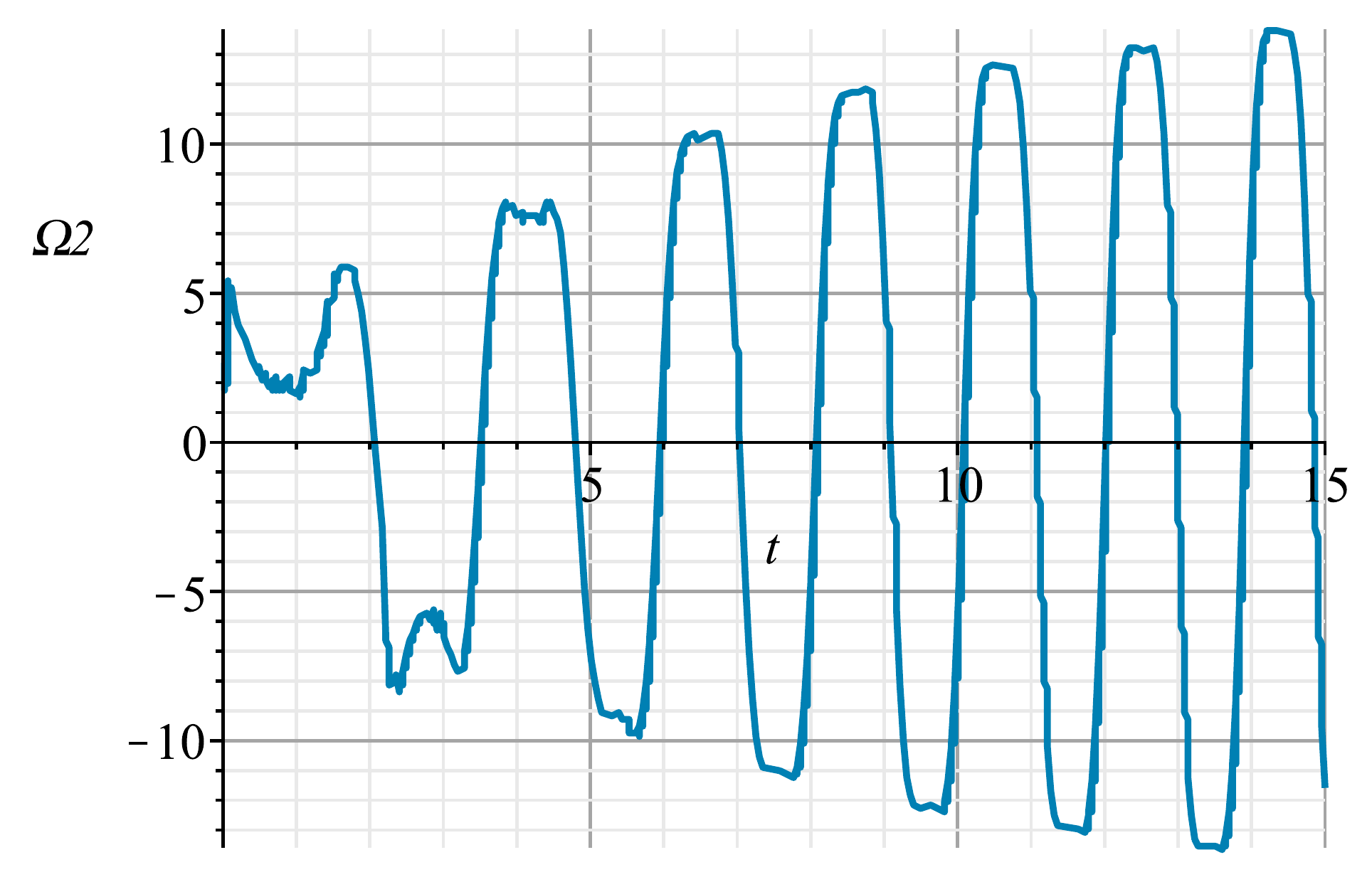}
}
\caption{Components $\Omega_i(t)$ of a sample path of the closed-loop system~\eqref{Euler2}, \eqref{controls2}.}
\label{fig2Omega}
\end{figure}

\begin{figure}
\center{
\includegraphics[scale=0.3]{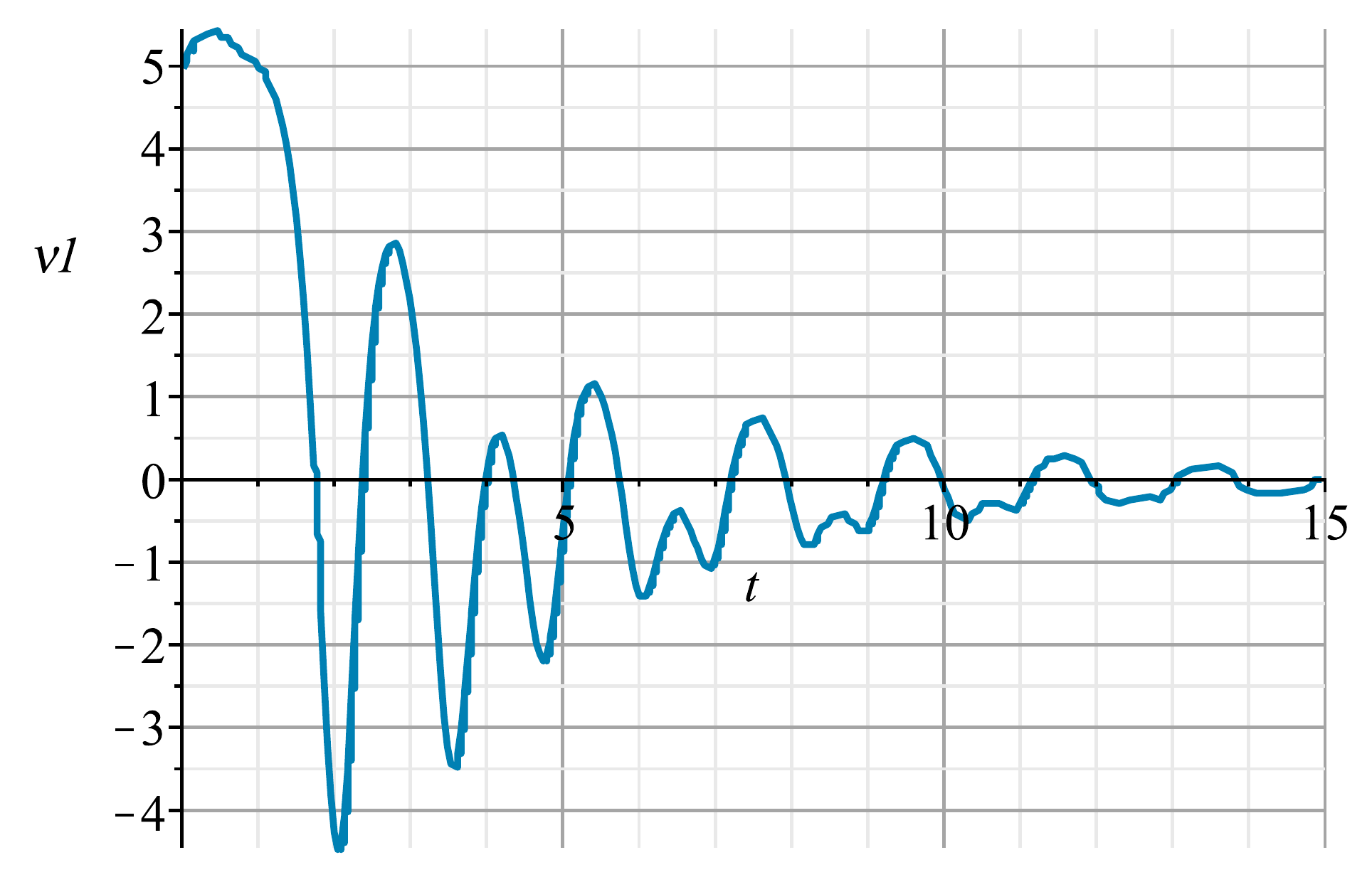}
\includegraphics[scale=0.3]{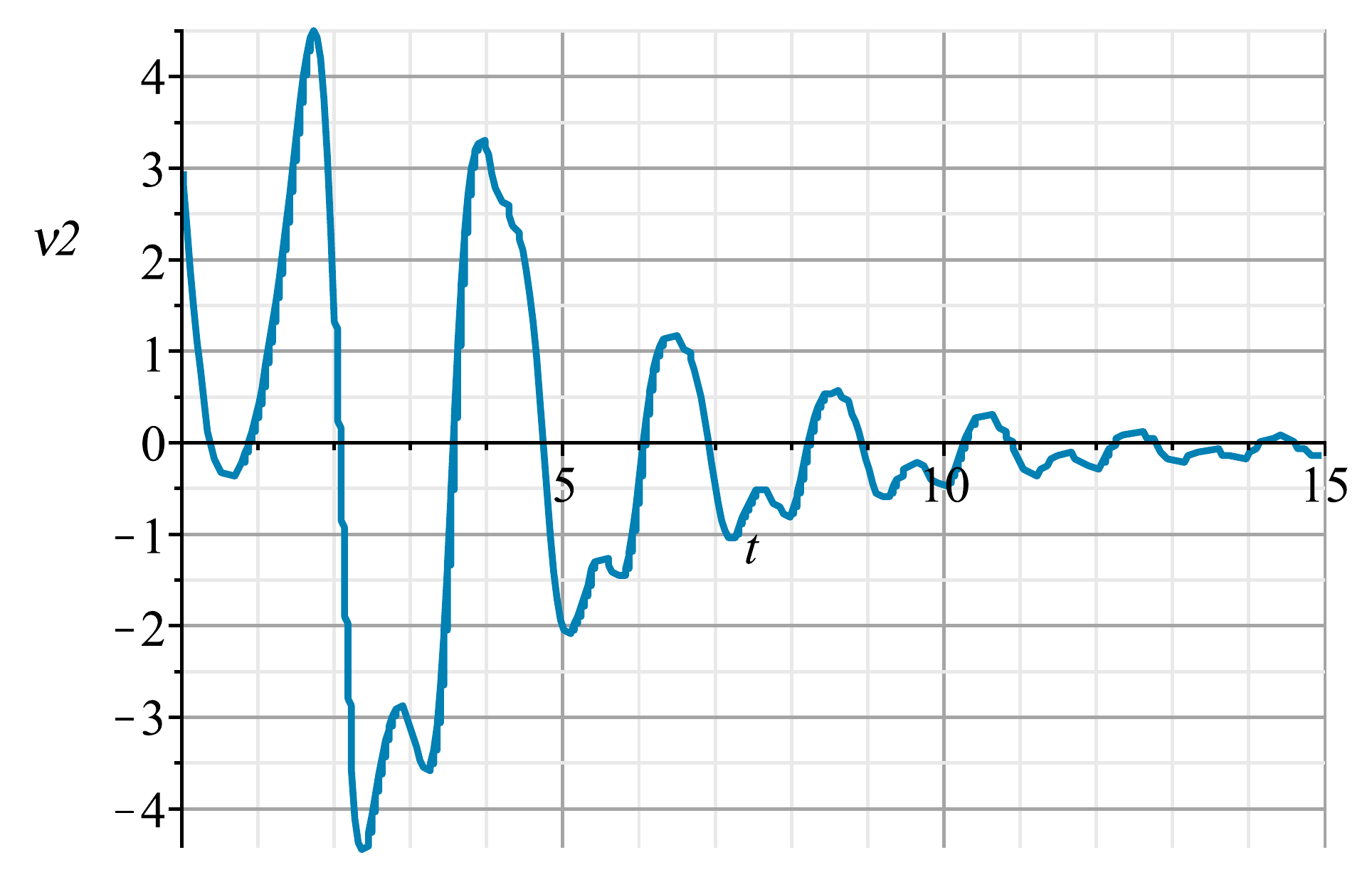}
\includegraphics[scale=0.3]{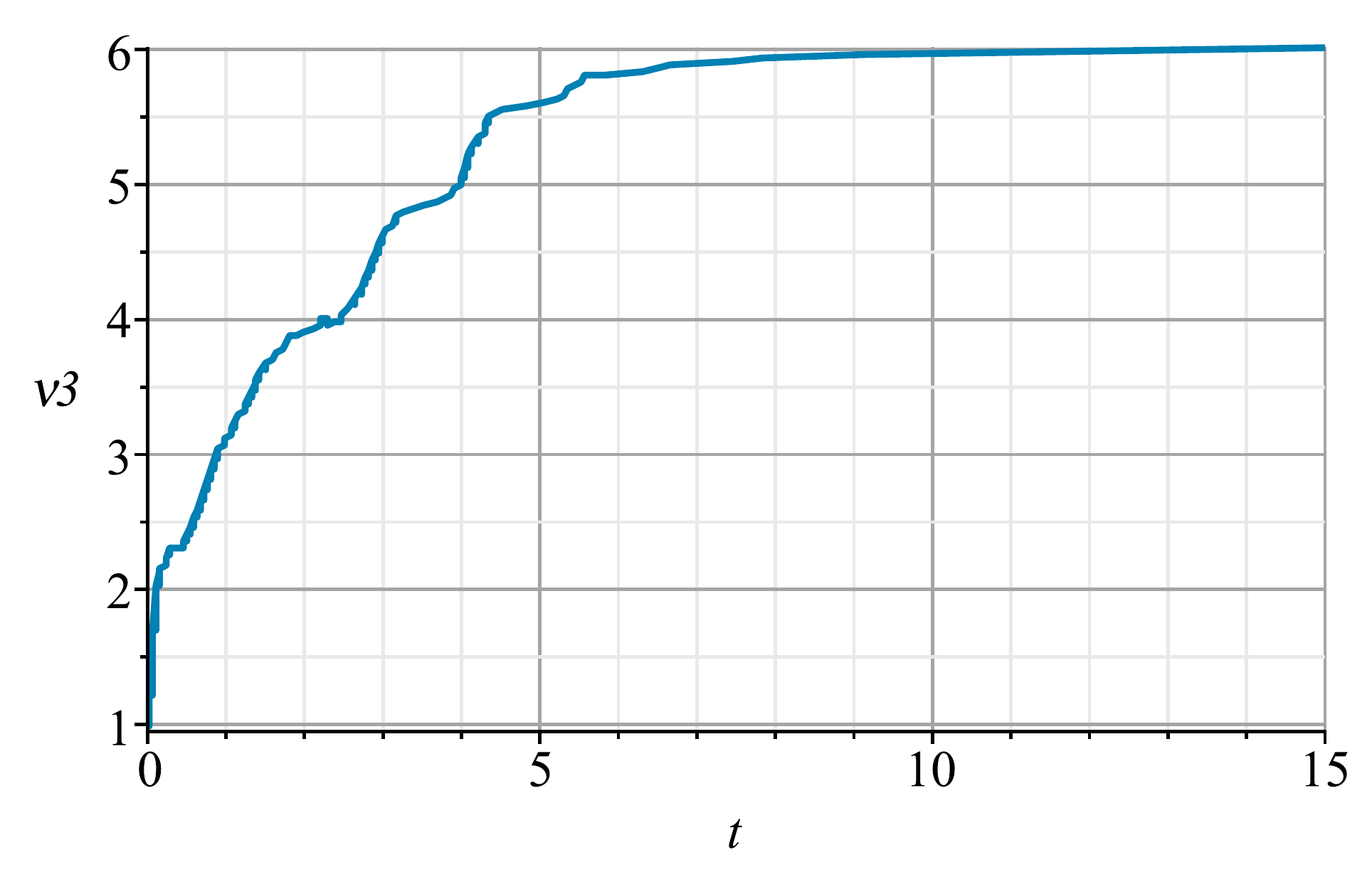}
}
\caption{Components $\nu_i(t)$ of a sample path of the closed-loop system~\eqref{Euler2}, \eqref{controls2}.}
\label{fig2nu}
\end{figure}

As in the previous example, we check the conditions of Theorem~1. We have: \small
$$\mathcal LV(x)=-\omega_1^2(h+\frac{|A_1-A_2|}{2}|\omega_3|)-\omega_2^2(h+\frac{|A_1-A_2|}{2}|\omega_3|)\leq0.$$\normalsize

By exploiting the integrals $W_1, W_2, W_3$ and the condition $\mathcal LV (x)\leq 0,$ we conclude about the boundedness of the solutions of system~\eqref{Euler2} with the feedback law~\eqref{controls2}.

The closed-loop system~\eqref{Euler2},~\eqref{controls2} on the set $M$ takes the form:
\begin{equation}
\begin{array}{lcl}
d\omega_{3}=0,\\
d\Omega_{1}=\Omega_2\omega_3dt,\\
d\Omega_{2}=\Omega_1\omega_3dt,\\
d\nu_{3}=0.
\end{array}
\end{equation}

Let us find the solutions of system (20):
\begin{equation}
\begin{array}{lcl}
\omega_{3}(t)=c_1\sqrt{\frac{(A-1-I_1)(A_2-I_2)}{I_1I_2}}\ne 0,\\
\nu_{3}(t)=c_2,\\
\Omega_{1}(t)=c_3 \cos(c_1t)+c_4 \sin(c_1t),\\
\Omega_{2}=\sqrt{\frac{I_2(A_2-I_2)}{I_1(A_1-I_1)}}c_4 \cos(c_1t)-c_3 \sin(c_1t).\\
\end{array}
\end{equation}

So, for arbitrary initial conditions from the set $M$, the expressions (21) define the solutions of system (20) in $M$, which proves the invariance of the set $M$.

The set $M_v$ has the form: $M_v=\{x\,|\,\omega_1=\omega_2=0\},$ i.e.
\begin{equation}x(t)\in M_v:
\begin{array}{lcl}
d\nu_1=\omega_3\nu_2dt,\\
d\nu_2=-\omega_3\nu_1dt,\\
\nu_2\nu_3=0,\\
\nu_1\nu_3=0,
d\omega_3=0,
d\nu_3=0.
\end{array}
\end{equation}

For the initial conditions sufficiently close to the equilibrium position (17), all the solutions of system \eqref{Euler2} with control \eqref{controls2} posess the property $\nu_1=\nu_2=0.$
Then the corresponding solution of system (22) is given by the relations $\nu_1(t) = 0, \nu_2(t) = 0$ together with (21).
The same solution satisfies system (20). This means that $M_v \backslash M$ does not contain any semi-trajectory of the considered closed-loop system, which means that the last condition of  Theorem~1  is satisfied.

Thus, the set $$M=\{(\omega_1,\omega_2,\omega_3,\Omega_1,\Omega_2,\nu_1,\nu_2,\nu_3)|\nu_i=\omega_i=0,\, i=1,2\}$$ is asymptotically stable  in probability for the closed-loop system \eqref{Euler2} with controls \eqref{controls2} by  Theorem~1.
This conclusion is also illustrated by  the results of  numerical simulations of the closed-loop system \eqref{Euler2}, \eqref{controls2} with the parameters $A_1=10, A_2=31, A_3=22, I_1=8, I_2=27, \varepsilon = 0.001.$

Figs.~\ref{fig2w}--\ref{fig2nu} show that the stabilized variables tend to zero as $t\to+\infty$.

\section{Conclusions}

In this paper, the idea of the Barbashin--Krasovskii theorem and LaSalle's invariance principle has been extended to characterize the asymptotic stability property of invariant sets of stochastic differential equations. It is shown that this result can be applied to the partial stabilization problem for nonlinear control systems with stochastic effects.
It should be emphasized that the control systems, considered in Sections~4 and~5, are not stabilizable in the classical sense with respect to all variables because of the presence of the geometric integral $\nu_1^2+\nu_2^2+\nu_3^2 = const$.
Thus, only partial stabilization is possible in the considered cases, and numerical simulations illustrate the efficiency of the proposed controllers.

\end{document}